\documentclass[12pt]{article}

\usepackage{amssymb}
\usepackage{amscd}
\usepackage[ansinew]{inputenc}
\usepackage[frenchb]{babel}

\input psfig.sty
\font\tencsc=cmcsc10

\def\ra{\rightarrow}

\def\lgl{\langle}
\def\rgl{\rangle}

\def\fl{\forall}
\def\ify{\infty}
\def\op{\oplus}
\def\ot{\otimes}

\def\part{\partial}
\def\sbs{\subset}
\def\sm{\simeq}
\def\ts{\times}
\def\ud{\underline}

\def\a{\alpha}
\def\d{\delta}
\def\D{\Delta}
\def\g{\gamma}
\def\G{\Gamma}
\def\lb{\lambda}

\def\s{\sigma}
\def\ve{\varepsilon}
\def\vp{\varphi}
\def\t{\theta}
\def\b{\beta}

\font\tenbb=msbm10
\font\sevenbb=msbm7
\font\fivebb=msbm5
\newfam\bbfam
\textfont\bbfam=\tenbb \scriptfont\bbfam=\sevenbb
\scriptscriptfont\bbfam=\fivebb
\def\bb{\fam\bbfam}

\def\Cb{{\bb C}}

\def\Rb{{\bb R}}

\def\Hc{{\cal H}}
\def\Lc{{\cal L}}
\def\Nc{{\cal N}}
\def\Oc{{\cal O}}
\def\Sc{{\cal S}}

\catcode`\@=11
\def\displaylinesno #1{\displ@y\halign{
\hbox to\displaywidth{$\@lign\hfil\displaystyle##\hfil$}&
\llap{$##$}\crcr#1\crcr}}

\def\ldisplaylinesno #1{\displ@y\halign{
\hbox to\displaywidth{$\@lign\hfil\displaystyle##\hfil$}&
\kern-\displaywidth\rlap{$##$}
\tabskip\displaywidth\crcr#1\crcr}}
\catcode`\@=12

\def\build#1_#2^#3{\mathrel{
\mathop{\kern 0pt#1}\limits_{#2}^{#3}}}

\def\Talpha#1{\vbox{\ialign{##\crcr
      $\alpha$\crcr\noalign{\kern2pt\nointerlineskip}
	   $\hfil\displaystyle{#1}\hfil$\crcr}}}

\def\Cb{{\bb C}}
\def\Fb{{\bb F}}

\def\Rb{{\bb R}}

\def\Cb{{\mathbb C}}

\def\Rb{{\mathbb R}}

\def\Hc{{\cal H}}
\def\Lc{{\cal L}}
\def\Nc{{\cal N}}
\def\Oc{{\cal O}}
\def\Sc{{\cal S}}

\def\Hc{{\cal H}}

\def\Lc{{\cal L}}

\def\Sc{{\cal S}}

\def\a{\alpha}
\def\b{\beta}
\def\d{\delta}
\def\lb{\lambda}
\def\g{\gamma}

\def\s{\sigma}
\def\t{\theta}
\def\ve{\varepsilon}
\def\vp{\varphi}

\def\D{\Delta}
\def\G{\Gamma}

\def\fl{\forall}
\def\ify{\infty}
\def\lgl{\langle}

\def\op{\oplus}
\def\ot{\otimes}

\def\part{\partial}
\def\rgl{\rangle}
\def\sbs{\subset}

\def\sm{\simeq}
\def\ts{\times}

\def\ra{\rightarrow}

\def\text{\hbox}

\def\Cb{{\mathbb  C}}

\def\build#1_#2^#3{\mathrel{
\mathop{\kern 0pt#1}\limits_{#2}^{#3}}}

\def\beq{\begin{equation}}
\def\eeq{\end{equation}}

\input{amssymb.sty}
\makeatletter
\def\nothing#1{}
\newdimen\earraycolsep
\renewcommand{\thetable}{\arabic{table}}
\renewcommand{\thefigure}{\arabic{figure}}
\renewcommand{\title}[1]{%
  \vspace*{50\p@}%
   {\parindent \z@ \raggedright \reset@font
     \bfseries #1\par
     \nobreak
     \vskip 36\p@
   }}
\def\author#1{{\pretolerance=10000 \raggedright \advance \leftskip by 2in
\noindent #1 \vskip 1pc}}
\def\affiliation#1{{\advance\leftskip by 1in \noindent #1 \vskip -1pc}}
\def\refnote#1{{$^{\hbox{\scriptsize #1}}$}}
\def\affnote#1{\llap{$^{\hbox{\scriptsize #1}}$}}

\renewcommand\section{\@startsection{section}{1}{\z@}{2pc \@plus
       1ex minus .2ex}{1pc \@plus .2ex}{\reset@font
       \normalsize\bfseries\noindent
       {\addtocounter{section}{1}}\Roman{section}\
       {\setcounter{subsection}{0}
       \setcounter{subsubsection}{0}\setcounter{equation}{0}} }}
\renewcommand\subsection{\@startsection{subsection}{2}{\z@}{1pc \@plus 1ex
     minus.2ex}{1pc \@plus .2ex}
     {\reset@font\normalsize\bfseries
     \noindent{\addtocounter{subsection}{1}}%
     {\setcounter{subsubsection}{0}}\Roman{section}.\Roman{subsection}\ }}
\renewcommand\subsubsection{\@startsection{subsubsection}{3}{\parindent}
         {1pc \@plus 1ex minus.2ex}{-0.5em}{\reset@font\normalsize\bfseries%
         {\addtocounter{subsubsection}{1}} \hspace*{.6cm}
         \Roman{section}.\Roman{subsection}.\Roman{subsubsection}
         \hspace*{-7mm}}}

\def\AmS{{\protect\the\textfont2%
         A\kern-.1667em\lower.5ex\hbox{M}\kern-.125emS}}

\def\p@LaTeX{{\family{times}\series{m}\shape{n}\selectfont
L\kern-.36em\raise.3ex\hbox{\scriptsize A}\kern-.15em
T\kern-.1667em\lower.7ex\hbox{E}\kern-.125emX}}

\newlength{\colwidth}

\def\@oddhead{\hfil}
\def\@evenhead{\hfil}
\def\@oddfoot{{\bfseries\hfil\thepage}}
\def\@evenfoot{{\bfseries\thepage\hfil}}
\def\fnum@figure{\footnotesize\raggedright{\bfseries \figurename~\thefigure.}}
\def\fnum@table{\normalsize\raggedright{\bfseries \tablename~\thetable.}}
\long\def\@makecaption#1#2{\vskip 10\p@ {#1 #2\par}}
\long\def\@makefntext#1{\setbox0=\hbox{$\m@th^{\@thefnmark}$}\noindent
\hangindent=\wd0 \box0 #1}
\newbox\@atbox
\long\def\atable#1#2#3{\begin{table}[tbp]\centering\footnotesize
\setbox\@atbox\hbox{#2}
\parbox{\wd\@atbox}{\caption{#1}}\par\smallskip
#2
\par\smallskip\parbox{\wd\@atbox}{\raggedright #3}
\end{table}}

\def\build#1_#2^#3{\mathrel{
\mathop{\kern 0pt#1}\limits_{#2}^{#3}}}

\def\@nbibitem#1{\noindent \hangindent=2pc \hangafter=1
\refstepcounter{enumi}\hbox to 2pc{\arabic{enumi}.\hfil}%
\immediate\write\@auxout{\string\bibcite{#1}{\arabic{enumi}}}}
\def\numbibliography{%
\section*{REFERENCES}%
\bgroup\footnotesize
\setcounter{enumi}{0}%
\def\newblock{\hskip .11em plus.33em minus.07em}%
\let\bibitem\@nbibitem}
\def\endnumbibliography{\par\egroup}
\makeatother
\def\ra{\rightarrow}

\font\tenbb=msbm10
\font\sevenbb=msbm7
\font\fivebb=msbm5
\newfam\bbfam
\textfont\bbfam=\tenbb \scriptfont\bbfam=\sevenbb
\scriptscriptfont\bbfam=\fivebb
\def\bb{\fam\bbfam}

\def\build#1_#2^#3{\mathrel{
\mathop{\kern 0pt#1}\limits_{#2}^{#3}}}

\begin{document}

\title{
\centerline{ SYM\'ETRIES GALOISIENNES ET}
\centerline{\, }
\centerline{ RENORMALISATION}
}

\author{\bf Alain CONNES\refnote{}}

\affiliation{\affnote{}Coll\`ege de France, 3, rue Ulm, 75005 PARIS\\
et \\ I.H.E.S., 35, route de Chartres, 91440 BURES-sur-YVETTE }

\begin{abstract}
\noindent Nous exposons notre travail en collaboration avec 
Dirk Kreimer sur les algèbres de Hopf et de Lie associées aux 
graphes de Feynman, et sur la signification conceptuelle
de la renormalisation perturbative à partir du problème 
de Riemann-Hilbert. Nous interprétons ensuite le groupe 
de renormalisation comme un groupe d'ambiguïté et 
montrons le rôle que ce 
groupe devrait jouer pour comprendre la composante  connexe du groupe des classes 
d'idèles de la théorie du corps de classe comme un groupe de Galois.
\end{abstract}

\section{Introduction}

\noindent La renormalisation est sans doute l'un des procédés les
plus élaborés pour obtenir des quantités numériques signifiantes à partir d'expressions mathématiques a-priori 
dépourvues de sens.
A ce titre elle est fascinante autant pour le physicien que 
pour le mathématicien. La 
profondeur de ses origines en théorie des champs et la précision 
avec laquelle elle est corroborée par l'expérience en font
l'un des joyaux de la physique théorique.
Pour le mathématicien épris de sens, mais non corseté par la rigueur,
les explications données jusqu'à présent butaient toujours sur
le sens conceptuel de la partie proprement calculatoire,
celle qui est utilisée par exemple en électrodynamique quantique
et ne tombe pas sous la coupe des "théories asymptotiquement libres"
auxquelles la théorie constructive peut prétendre avoir donné un statut
mathématique satisfaisant.
Cet état de fait a changé récemment et cet exposé se propose de
donner la signification conceptuelle 
des calculs effectués par les physiciens dans la théorie de
la renormalisation grâce
à mon travail sur 
la renormalisation en 
collaboration avec Dirk Kreimer et la relation que nous avons établie entre 
renormalisation et problème
de Riemann-Hilbert.

\noindent Le résultat clef est l'identité entre le procédé récursif utilisé par les physiciens et les formules 
mathématiques qui résolvent un lacet à valeurs dans un groupe pronilpotent $G$ en un produit d'un lacet 
holomorphe par un lacet anti-holomorphe. La signification géométrique de cette décomposition (de 
Riemann-Hilbert, Birkhoff ou Wiener-Hopf) provient directement de la théorie des fibrés holomorphes
de groupe structural $G$ sur la sphère de Riemann $S^2$. Dans la renormalisation perturbative, les points de la 
sphère $S^2$ sont les dimensions complexes parmi lesquelles la dimension $d$ de l'espace-temps est un point 
privilégié. Le problème étant que dans les théories physiquement intéressantes les quantités à calculer 
conspirent pour diverger précisément au point $d$.  On peut organiser ces quantités comme le développement de 
Taylor d'un difféomorphisme $g \in G$ et donner un sens à $g=g(z)$
en  remplaçant dans les formules la dimension $d$ par une valeur complexe $z \neq d$. Le procédé de 
renormalisation acquiert alors la signification suivante : la valeur cherchée $g \in G$ n'est autre que la 
valeur $g_+(d)$ en $d$ 
de la partie holomorphe de la décomposition de Riemann-Hilbert $ g(z)= g_-^{-1}\,(z) g_+(z)$ du lacet $ g(z)$. 
La nature exacte du groupe $G$ impliqué dans la renormalisation a été clarifiée par les étapes essentielles  
suivantes. 
La première est la découverte due à Dirk Kreimer de la structure d'algèbre de Hopf secrètement présente dans 
les formules récursives de Bogoliubov Parasiuk Hepp et Zimmermann.

\noindent La seconde qui est le point de départ de notre collaboration est
la similitude entre l'algèbre de Hopf des arbres enracinés de Dirk et une algèbre de Hopf que j'avais 
introduite avec Henri Moscovici pour
organiser les calculs très complexes de géométrie noncommutative.
Ceci nous a conduit avec Dirk à définir une algèbre de Hopf directement
en termes de graphes de Feynman et à lui appliquer le théorème de 
Milnor-Moore pour en déduire une algèbre de Lie et un groupe de Lie pronilpotent $G$, analogue
du groupe des difféomorphismes formels.

\noindent Enfin la troisième étape cruciale est la construction d'une
action du groupe $G$ sur les constantes de couplage de la théorie physique. 
Ceci permet de relever le groupe de renormalisation comme un sous-groupe
à un paramètre du groupe $G$ et de montrer directement que les développements
polaires des divergences sont entièrement déterminés  par leurs résidus.

\noindent Le problème de Riemann-Hilbert joue un rôle clef dans la théorie de Galois différentielle, il est 
donc naturel d'interpréter en termes Galoisiens l'ambiguïté que le groupe de renormalisation introduit dans les 
théories physiques.  La dernière section contient l'esquisse d'une telle interprétation. 

\noindent Nous commencerons cette section par une introduction
 très élémentaire à la théorie de Galois pour les équations algébriques,
en passant par un beau problème de géométrie plane.

\noindent Nous montrerons ensuite le rôle que le 
groupe de renormalisation
devrait jouer pour comprendre la composante  connexe du groupe des classes 
d'idèles de la théorie du corps de classe comme un groupe de Galois.
Cette idée s'appuie à la fois sur l'analogie entre la théorie des facteurs et 
la théorie de Brauer pour un corps local et sur la présence
implicite en théorie des champs d'un "corps de constantes" plus élaboré que le corps $\Cb$
des nombres complexes. En fait les calculs des physiciens regorgent d'exemples de "constantes"
telles les constantes de couplage $g$ des interactions (électromagnétiques, 
faibles et fortes)
 qui n'ont de "constantes"
que le nom. Elles dépendent en réalité du niveau d'énergie $\mu$ 
auquel les expériences sont réalisées et sont des fonctions $g( \mu)$, de sorte que les physiciens des hautes 
énergies
étendent implicitement le "corps des constantes" avec lequel
ils travaillent, passant du corps $\Cb$ des scalaires 
à un corps de fonctions $g( \mu)$. Le groupe d'automorphismes
de ce corps engendré par $ \mu \partial 
/ \partial \mu$
 est le groupe d'ambiguïté de la 
théorie physique.

\section { Renormalisation, position du problème}

\noindent La motivation physique de la renormalisation est très claire et remonte 
aux travaux de 
Green au dix-neuvième siècle sur l'hydrodynamique. Pour prendre un exemple simple 
\footnote{{\it voir le 
cours de théorie des champs de Sidney Coleman}} si l'on calcule l'accélération 
initiale d'une balle
 de ping-pong plongée à quelques mètrès sous l'eau, l'on obtient en appliquant la 
loi de Newton
$F=ma$ et la poussée d'Archimède $F= (M-m)g$, où $m$ est la masse inerte, et $M$ 
la masse
d'eau occupée, une accélération initiale de l'ordre de $11g$! \footnote{{\it La 
balle pèse $m=2,7$ grammes
 et a un diamètre de $4$ cm de sorte que $M=33,5$ grammes}}
 En réalité, si l'on réalise l'expérience,
l'accélération est de l'ordre de $2g$. En fait la présence du fluide autour de la 
balle oblige à 
corriger la valeur $m$ de la masse inerte dans la loi de Newton et à la remplacer 
par une 
"masse effective" qui en l'occurrence vaut $m + $ ${1} \over{ 2}$ $M$.
 Dans cet exemple, l'on peut bien sur déterminer la masse $m$ en pesant la balle 
de ping-pong hors de l'eau,
mais il n'en va pas de même pour un electron dans le champ electromagnétique, dont 
il est impossible de l'extraire.
De plus le calcul montre que, pour une particule ponctuelle comme le demande la 
relativité, la correction qui valait
${1} \over{ 2}$ $M$ ci-dessus est infinie.

\noindent Vers 1947 les physiciens ont réussi à utiliser la distinction entre les 
deux masses qui apparaissent ci-dessus et
plus généralement le concept de quantité physique "effective" pour éliminer les 
quantités infinies qui apparaissent
en théorie des champs quantiques (voir \cite{dresden} pour un aperçu historique).

\noindent Une th\'eorie des champs en $d$ dimensions est donn\'ee par une 
fonctionnelle d'action classique
 $$
S \, (A) =  \int \Lc \, (A) \,
d^d x \leqno(1) 
$$ 
o\`u $A$ d\'esigne un champ classique et le 
Lagrangien est de la forme, 
$$ \Lc \, (A) =  (\part A)^2 / 2 -
\frac{m^2}{2} \, A^2 + \Lc_{\rm int} (A) \leqno(2)
$$ 
o\`u
$\Lc_{\rm int} (A)$ est un polynome en $A$.

\smallskip

\noindent On peut d\'ecrire la th\'eorie par les fonctions de Green,
$$
G_N (x_1 , \ldots , x_N) = \lgl \, 0 \, \vert T \, \phi (x_1) \ldots
\phi (x_N) \vert \, 0 \, \rgl \leqno(3)
$$
o\`u le symbole $T$ signifie que les champs quantiques $\phi (x_j)$'s sont écrits \`a temps 
croissant de droite \`a gauche.

\smallskip

\noindent L'amplitude de probabilit\'e d'une configuration classique $A$ est 
donn\'ee par,
$$
e^{i \, \frac{S(A)}{\hbar}} \leqno(4)
$$
et si l'on pouvait ignorer les probl\`emes de renormalisation, l'on pourrait 
calculer les fonctions de Green gr\^ace \`a la formule
$$
G_N (x_1 , \ldots , x_N) = \Nc \int e^{i \, \frac{S(A)}{\hbar}} \ A (x_1)
\ldots A (x_N) \, [dA] \leqno(5)
$$
o\`u $\Nc$ est un facteur de normalisation requis par,  
$$
\lgl \, 0 \mid 0 \, \rgl = 1 \, . \leqno(6)
$$

\smallskip

\noindent  L'on pourrait alors calculer 
l'int\'egrale fonctionnelle (5) en th\'eorie des perturbations, en traitant le 
terme $\Lc_{\rm int}$ de (2) 
comme une perturbation, le Lagrangien libre étant,
$$
\Lc_0 (A) = (\part A)^2 / 2 - \frac{m^2}{2} \, A^2 \, ,
\leqno(7)
$$
de sorte que,
$$
S (A) = S_0 (A) + S_{\rm int} (A) \leqno(8)
$$
o\`u l'action libre $S_0$ d\'efinit une mesure Gaussienne $\exp \, (i
\, S_0 (A)) \, [dA] = d \mu$.

\smallskip

\noindent On obtient alors le d\'eveloppement perturbatif des fonctions de Green 
sous la forme,
$$
\ldisplaylinesno{
G_N (x_1 , \ldots , x_N) = \left( \sum_{n=0}^{\ify} \, i^n / n! \int
A (x_1) \ldots A (x_N) \, (S_{\rm int} (A))^n \, d \mu
\right) \cr
\quad \left( \sum_{n=0}^{\ify} \, i^n / n! \int S_{\rm int} (A)^n \, d
\mu \right)^{-1} \, . &(9) \cr
}
$$
Les termes de ce d\'eveloppement s'obtiennent en int\'egrant par partie sous la 
Gaussienne.
Cela engendre un grand nombre de termes $U(\G)$, où les paramètrès $\G$  sont 
les graphes de Feynman $\G$, i.e. des graphes dont les sommets correspondent 
aux termes du Lagrangien de la th\'eorie. 
\smallskip
En règle générale les valeurs des termes $U(\G)$ sont
donn\'ees par des int\'egrales divergentes. Les divergences les plus importantes 
sont causées
par la présence dans le domaine d'intégration de moments de taille arbitrairement 
grande.
La technique de renormalisation consiste d'abord à "régulariser" ces int\'egrales 
divergentes
par exemple en introduisant un paramètre de "cutoff" $\Lambda$ et en se 
restreignant à la portion
correspondante du domaine  d'intégration. Les intégrales sont alors finies, mais 
continuent bien entendu
à diverger quand $\Lambda \rightarrow \infty$. On établit ensuite une dépendance 
entre les termes
du Lagrangien et $\Lambda$ pour que les choses s'arrangent et que les résultats 
ayant un sens physique
deviennent finis! 
Dans le cas 
particulier des théories asymptotiquement libres, la forme explicite de la dépendance
 entre les constantes nues et le paramètre de régularisation $\Lambda$ a permis dans des cas 
très importants (\cite{GK},\cite{R}) de mener à bien le programme de la théorie constructive
des champs (\cite{GJ}). 

\noindent Décrivons maintenant en détail la technique de renormalisation perturbative.
Pour faire les choses systématiquement, on rajoute un "contreterme" 
$C(\G)$ au Lagrangien de départ ${\cal L}$,
 chaque fois que l'on rencontre un diagramme divergent $\G$, dans le but 
d'annuler la divergence correspondante. Pour les théories "renormalisables", les 
contre-termes
dont on a besoin sont tous déja des termes du Lagrangien ${\cal L}$ et ces 
contorsions peuvent
s'interpreter à partir de l'inobservabilité des quantités numériques qui 
apparaissent dans ${\cal L}$,
par opposition aux quantités physiques qui, elles, doivent rester finies.

\noindent La principale complication dans cette procédure vient de 
l'existence de nombreux graphes $\G$ pour lesquels les divergences de $U(\G)$ ne
sont pas locales. La raison étant que ces graphes possèdent déja des sous-graphes
dont les divergences doivent être prises en compte avant d'aller plus loin.
La méthode combinatoire précise, due à Bogoliubov-Parasiuk-Hepp et Zimmermann 
(\cite{bphz})
consiste d'abord à "préparer" le graphe $\G$ en remplaçant $U(\G)$ par 
l'expression
formelle,
$$
\overline{R}( \G) = U(\G) + \sum_{\g \sbs \G} C(\g) U( \G /
\g) \leqno (10)
$$
où $\g$ varie parmis tous les sous-graphes divergents.
On montre alors que le calcul des divergences du graphe "préparé" ne donne que des 
expressions
locales, qui pour les théories renormalisables se trouvent déja dans le Lagrangien 
 ${\cal L}$.

\section  { L'alg\`ebre de Hopf des graphes de Feynman}

\noindent Dirk Kreimer a eu l'idée remarquable en 97 (\cite{K1}) d'utiliser la 
formule (10) pour définir 
le coproduit d'une algèbre de Hopf.

\noindent  En tant qu'alg\`ebre $\Hc$ est l'alg\`ebre commutative libre 
engendr\'ee par les graphes "une particule irréductibles" (1PI)\footnote{{\it  Un 
graphe de Feynman  $\G$ est "une particule irreductible"
(1PI) si il est connexe et le reste apres avoir enlev\'e n'importe laquelle de ses 
faces}}

\noindent Elle admet ainsi une base index\'ee par les graphes $\G$ unions 
disjointes de graphes 1PI.
$$
\G = \bigcup_{j=1}^{n} \G_j \, . \leqno (11)
$$

\noindent Le produit dans $\Hc$ est donn\'e par l'union disjointe,
$$
\G \cdot \G' = \G \cup \G' \, . \leqno (12)
$$
Pour d\'efinir le coproduit,
$$
\D : \Hc \ra \Hc \ot \Hc \leqno (13)
$$
il suffit de le donner sur les graphes 1PI, on a
$$
\Delta \G = \G \ot 1 + 1 \ot \G + \sum_{\g \sbs \G} \g_{(i)} \ot \G /
\g_{(i)} \leqno (14)
$$

\smallskip

\noindent Ici $\g$ est un sous-ensemble (non vide et de compl\'ementaire non-vide) 
$\g \sbs \G^{(1)}$ 
de l'ensemble $ \G^{(1)}$ des faces internes de $\G$ dont les composantes connexes 
$\g'$ 
v\'erifient des conditions d'admissibilit\'e d\'etaill\'ees dans la r\'ef\'erence 
\cite{ck3}.
\smallskip

\noindent Le coproduit $\D$ defini par (14) sur les graphes 1PI se prolonge de 
mani\`ere unique en un 
homomorphisme de $\Hc$ dans $\Hc \ot \Hc$. On a alors (\cite{K1},\cite{ck3})

\medskip

\noindent {\tencsc Théorème } {\it Le couple  $(\Hc , \D)$ est une alg\`ebre de
Hopf.}

\section  { L'algèbre de Lie des graphes, le groupe $G$ et sa structure.}
\medskip
\noindent  J'avais à la mème époque, dans les calculs de Géométrie
Noncommutative de l'indice transverse pour les feuilletages, montré, avec Henri 
Moscovici, (\cite{cm1}) que 
la compléxité extrême de ces calculs conduisait à introduire une algèbre de Hopf  
$\Hc_{cm}$, qui n'est 
ni commutative ni cocommutative mais est intimement reliée au groupe des 
difféomorphismes,
dont l'algèbre de Lie apparait en appliquant le théorème de Milnor-Moore à une 
sous-algèbre commutative.

\noindent Aprés l'exposé de Dirk à l'IHES en février 98, nous avons tous les deux 
été intrigués par la 
similarité apparente entre ces deux algèbres de Hopf et notre 
 collaboration a commencé par l'application du théorème de Milnor-Moore à 
l'algèbre de Hopf  $\Hc$.
Le th\'eor\`eme de Milnor-Moore montre qu'elle est duale de l'alg\`ebre 
enveloppante 
d'une alg\`ebre de Lie gradu\'ee $\ud G$ dont une base est donn\'ee par les 
graphes 1-particule irreductibles. Le crochet de Lie de deux graphes
 est obtenu par insertion d'un graphe dans l'autre. Le groupe de Lie correspondant 
$G$ est le groupe des caractères de $\Hc$. 

\noindent Nous avons ensuite analysé le groupe $G$ et montré qu'il est produit 
semi-direct d'un groupe ab\'elien par un groupe 
très reli\'e au groupe des diff\'eomorphismes des constantes de couplage sans 
dimension de la th\'eorie des champs
(voir section VII).
\smallskip
 
\noindent L'alg\`ebre de Hopf $\Hc$ admet plusieurs graduations naturelles.
Il suffit de donner le degr\'e des graphes 1PI puis de poser en g\'en\'eral,
$$
\deg \, (\G_1 \ldots \G_e) = \sum \deg \, (\G_j) \ , \quad \deg \, (1) = 0
\leqno (15)
$$
On doit v\'erifier que,
$$
\deg \, (\g) + \deg \, (\G / \g) = \deg \, (\G) \leqno (16)
$$
pour tout sous-graphe admissible $\g$.

\smallskip

\noindent Les deux graduations les plus naturelles sont
$$
I \, (\G) = \hbox{nombre de faces internes} \ \G \leqno (17)
$$
et
$$
v \, (\G) = V \, (\G) - 1 = \hbox{nombre de sommets} \ \G - 1 \, .
\leqno (18)
$$
On a aussi la combinaison importante
$$
L = I - v = I - V + 1 \leqno (19)
$$
qui est le nombre de boucles du graphe.

\noindent Soit $G$ un graphe 1PI avec $n$ faces externes index\'ees par
$i \in \{ 1 , \ldots , n \}$, on sp\'ecifie sa structure externe en donnant
une distribution $\s$ definie sur un espace convenable de fonctions test $\Sc$
sur $$ \left\{ (p_i)_{i=1 , \ldots , n} \ ; \
\sum \, p_i = 0 \right\} = E_G \, . \leqno (20) $$ Ainsi $\sigma$
est une forme lin\'eaire continue, $$ \sigma : \Sc \, (E) \ra \Cb \, .
\leqno (21) $$ A un graphe $\G$ de structure externe $\sigma$ correspond un 
\'el\'ement de $\Hc$ 
et on a $$ \d_{(\G , \lb_1
\sigma_1 + \lb_2 \sigma_2)} = \lb_1 \, \d_{(\G , \sigma_1)} +
\lb_2 \, \d_{(\G , \sigma_2)} \, . \leqno (22) $$ 
\smallskip

\noindent Nous appliquons alors le th\'eor\`eme de Milnor-Moore \`a l'alg\`ebre de 
Hopf bigradu\'ee
 $\Hc$. 

\smallskip

\noindent Ce th\'eor\`eme donne une structure d'alg\`ebre de lie sur,
 $$ \bigoplus_{\G} \ \Sc \, (E_{\G}) = L \leqno (23)
$$ o\`u pour chaque graphe 1PI $\G$, on d\'efinit $\Sc \, (E_{\G})$ comme
 dans (20).
Soit $X \in L$ et soit $Z_X$ la forme lin\'eaire sur $\Hc$ donn\'ee,
sur les monomes $\G$, par $$ \lgl \G , Z_X \rgl = 0 \leqno (24) $$
sauf si $\G$ est connexe et 1PI, et dans ce cas par, $$ \lgl \G , Z_X
\rgl =\lgl \sigma_{\G} , X_{\G} \rgl \leqno (25) $$ o\`u
$\sigma_{\G}$ est la distribution qui donne la structure externe de
$\G$ et $X_{\G}$ la composante correspondante de $X$. Par
construction $Z_X$ est un caractère infinit\'esimal de $\Hc$ 
ainsi que les commutateurs, $$ [Z_{X_1} , Z_{X_2}] =
Z_{X_1} \, Z_{X_2} - Z_{X_2} \, Z_{X_1}. \leqno (26) $$
Le produit étant obtenu par transposition du coproduit de $\Hc$,
i.e. par $$ \lgl Z_1 \, Z_2 , \G \rgl = \lgl Z_1 \ot Z_2 , \D \, \G
\rgl \, . \leqno (27) $$

\noindent  Soient $\G_j$, $j = 1,2$ des graphes 1PI et $\vp_j
\in \Sc \, (E_{\G_j})$ les fonctions test correspondantes.
\smallskip

\noindent Pour $i \in \{ 0,1 \}$, soit $n_i \, (\G_1 , \G_2 ; \G)$ le
nombre de sous-graphes de $\G$ isomorphes \`a $\G_1$ et tels que
$$
\G / \G_1 (i) \sm \G_2 \, . \leqno (28)
$$
Soit $(\G , \vp)$ l'élément de $L$ associ\'e \`a $\vp \in \Sc \,
(E_{\G})$, le crochet de Lie de $(\G_1 , \vp_1)$ et $(\G_2 , \vp_2)$ est donn\'e 
par,
$$
\sum_{\G , i} \, \sigma_i \, (\vp_1) \ n_i \, (\G_1 , \G_2 ; \G) \, (\G ,
\vp_2)
- \sigma_i \, (\vp_2) \ n_i \, (\G_2 , \G_1 ; \G) \, (\G , \vp_1) \, . \leqno
(29)
$$

\medskip

\noindent {\tencsc Théorème (\cite{ck3})} {\it L'alg\`ebre de Lie $L$ est produit 
semi-direct
d'une alg\`ebre de Lie Abelienne $L_0$ par $L_c$ o\`u $L_c$ admet une base 
canonique
 ind\'ex\'ee par les graphes $\G^{(i)}$ avec
$$
[\G , \G'] = \sum_v \, \G \circ_v \G' - \sum_{v'} \, \G' \circ_{v'} \G
$$
o\`u $\G \circ_v \G'$ est obtenu en greffant $\G'$ sur $\G$ en $v$.}

\medskip

\section  { Renormalisation et problème de Riemann-Hilbert}

\noindent  Le problème de Riemann-Hilbert vient du 
$21^{eme}$ problème de Hilbert qu'il formulait ainsi; 

 \begin{itemize}
\item[] ``Montrer qu'il existe toujours une équation différentielle Fuchsienne 
linéaire
de singularités et  monodromies données.''
\end{itemize}

\noindent Sous cette forme il admet une réponse positive  due à Plemelj
et Birkhoff (cf.~\cite{B} pour un exposé détaillé). Quand on le reformule pour les 
 systèmes linéaires de la  forme, $$ y'(z) = A(z)
\, y(z) \ , \ A(z) = \sum_{\a \in S} \ \frac{A_{ \a }}{z - \a} \,
, \leqno (30) $$ où $S$ est l'ensemble fini donné des
singularités, $\ify \not\in S$, et les $A_{\a}$ sont des  matrices complexes
telles que $$ \sum \ A_{\a} = 0 \leqno (31) $$ pour éviter les 
singularités à $\ify$, la réponse n'est pas toujours positive 
\cite{Bo}, mais la solution existe quand les  matrices de monodromie
$M_{\a}$  sont suffisamment proches de 1. On peut alors 
l'écrire explicitement sous la forme d'une série de 
 polylogarithmes \cite{LD}.

\noindent Une autre formulation du problème de Riemann-Hilbert,
intimement reliée à la classification des fibrés vectoriels holomorphes 
sur la sphère de Riemann $P_1 (\Cb)$, est en termes de la décomposition de
Birkhoff $$ \g \, (z) = \g_- (z)^{-1} \, \g_+ (z)
\qquad z \in C \leqno (32) $$ où $C \sbs P_1 (\Cb)$ désigne une courbe simple, 
$C_-$ la composante connexe du complément de $C$
contenant $\ify \not\in C$ et $C_+$ la composante bornée. 

$$
\hbox{ 
\psfig{figure=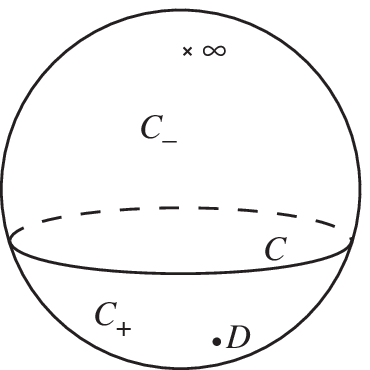} 
} 
\label{fig6}
$$

\centerline{Figure 1}

\medskip

\noindent Les trois lacets
$\g$ et $\g_{\pm}$ sont à valeurs dans ${\rm GL}_n (\Cb)$,
$$ \g \, (z) \in G = {\rm GL}_n (\Cb) \qquad \fl \, z \in \Cb
\leqno (33) $$ et $\g_{\pm}$ sont les valeurs au bord d'applications  holomorphes
$$ \g_{\pm} : C_{\pm} \ra {\rm GL}_n (\Cb) \, .
\leqno (34) $$
La condition $\g_- (\ify) = 1$
assure l'unicité de la décomposition (32) si elle existe.

\smallskip

\noindent L'existence de la décomposition de Birkhoff (32) est
équivalente à l'annulation,
$$
c_1 \, (L_j) = 0 \leqno (35)
$$
des nombres de Chern $n_j = c_1 \, (L_j)$ des fibrés en droites holomorphes
de la décomposition de Birkhoff-Grothendieck,
$$
E = \op \, L_j \leqno (36)
$$
où $E$ est le fibré vectoriel holomorphe sur $P_1 (\Cb)$ associé à
$\g$, i.e. d'espace total:
$$
(C_+ \ts \Cb^n)\cup_{\g} (C_- \ts \Cb^n) \, . \leqno (37)
$$
La discussion ci-dessus pour $G = {\rm GL}_n (\Cb)$ s'étend aux groupes de Lie 
complexes 
arbitraires. 
\smallskip

\noindent Quand $G$ est un groupe de Lie complexe nilpotent et simplement connexe
l'existence (et l'unicité) de la décomposition de  Birkhoff  (32)
est vraie pour tout $\g$. Quand le lacet $\g : C \ra G$ se prolonge en un lacet
holomorphe: $C_+ \ra G$, la décomposition de  Birkhoff est donnée par
$\g_+ = \g$, $\g_- = 1$. En général, pour $z \in C_+$ l'évaluation,
$$
\g \ra \g_+ (z) \in G \leqno (38)
$$
donne un  principe naturel pour extraire une valeur finie à partir de 
l'expression singulière $\g (z)$. Cette extraction de partie finie est
une division par la partie polaire 
pour un lacet méromorphe
$\g$ en prenant pour
 $C$ un cercle  infinitésimal centré en $z_0$.

\noindent Soit $G$ un groupe de Lie complexe pro-nilpotent, ${\mathcal H}$ son
algèbre de  Hopf de coordonnées (graduée). Rappelons la traduction entre langages 
algébriques et géométriques,
en désignant par  $\mathcal R$ l'anneau des fonctions méromorphes, ${\mathcal 
R}_-$ le sous anneau des polynomes en $(z-z_0)^{-1}$ et ${\mathcal R}_+$ celui des 
fonctions régulières en $z_0$,

 {\small
\begin{equation}
\begin{array}{rcl}
\underline{\mbox{{ Homomorphismes}  de ${\mathcal H}\to {\mathcal
R}$}} & \mid & \underline{\mbox{{ Lacets}  de $C$ à valeurs dans $G$}}\\
 & \mid & \\
 \phi({\mathcal
H})\subset\;{\mathcal R}_- & \mid & \mbox{$\gamma$ se prolonge en 
une application 
holomorphe de $\Cb\backslash\{z_0\}\to G$}\\
 & \mid &
\mbox{avec $\gamma(\infty)=1$.}\\ & \mid & \\ \phi({\mathcal
H})\subset\;{\mathcal R}_+ & \mid & \mbox{$\gamma$ se prolonge en 
une application 
holomorphe  définie en $z=z_0$.}\\ & \mid & \\ \phi=
\phi_1\star\phi_2 & \mid & \gamma(z)=\gamma_1(z)\gamma_2(z),\;
\forall\; z\in {\bf C}.\\ & \mid & \\
 \phi\circ S & \mid & z\to\gamma(z)^{-1}.
\end{array}
\end{equation}}

\noindent Pour $X \in{\mathcal H}$ notons le coproduit sous la forme
 $$ \Delta(X)= X \otimes 1 + 1 \otimes X + \sum X^\prime
\otimes X^{\prime\prime} $$
 La décomposition de Birkhoff
d'un lacet s'obtient de manière récursive grâce au théorème suivant,

\medskip

\noindent  {\tencsc Théorème (\cite{ck3})} {\it Soit $\phi: {\mathcal H}\to
{\mathcal R}$, un homomorphisme d'algèbres. La  décomposition de Birkhoff du lacet 
correspondant est donnée de manière récursive par les égalités,
$$\phi_-(X)=-T\left(\phi(X)+\sum\phi_-(X^\prime)
\phi(X^{\prime\prime}) \right)$$
$$\phi_+(X)=\phi(X)+\phi_-(X)+\sum\phi_-(X^\prime)
\phi(X^{\prime\prime}).$$}

\medskip

\noindent Ici $T$ désigne la projection sur ${\mathcal R}_-$ parallèlement à 
${\mathcal R}_+$.

\noindent La clef de notre travail avec Dirk Kreimer réside dans l'identité entre
ces formules et celles qui gouvernent la combinatoire des calculs de graphes. Nous 
avons déja vu
la formule qui définit la préparation d'un graphe,

$$
\overline{R}( \G) = U(\G) + \sum_{\g \sbs \G} C(\g) U( \G /
\g) \leqno (39)
$$
Celle qui donne le contreterme $C(\G)$ est alors,

$$
C(\G) =  -T(\overline{R}( \G)) = -T\left(U(\G) + \sum_{\g \sbs \G} C(\g) U( \G /
\g)\right) \leqno (40)
$$
et celle qui donne la valeur renormalisée du graphe est,

$$
R(\G) =  \overline{R}( \G) +C(\G) =U(\G)  +C(\G) + \sum_{\g \sbs \G} C(\g) U( \G /
\g)  \leqno (41)
$$
Il est alors clair en posant $\phi=U$, $\phi_-=C$, $\phi_+=R$ que ces équations 
sont identiques à celles du théorème
donnant la construction récursive de la décomposition de Birkhoff.

\medskip
\noindent Décrivons plus en détails ce résultat. Etant donn\'ee une th\'eorie 
renormalisable en dimension $D$
la th\'eorie nonrenormalis\'ee donne en utilisant la r\'egularisation 
dimensionnelle un lacet 
$\g$ d'éléments du groupe $G$ associ\'e \`a la th\'eorie dans la section IV. 
Le paramètre $z$ du lacet $\g \, (z)$ est une variable complexe
et $\g \, (z)$ est m\'eromorphe dans un voisinage de $D$.
Notre r\'esultat principal est que la th\'eorie renormalis\'ee est donn\'ee par 
l'\'evaluation \`a $z=D$
de la partie nonsingulière $\g_+$ de la d\'ecomposition de Birkhoff ,
$$
\g \, (z) = \g_- \, (z)^{-1} \ \g_+ \, (z) \leqno (42)
$$
de $\g$.
\smallskip

\noindent Les r\`egles de Feynman et la r\'egularisation dimensionnelle associent 
un nombre,
$$
U_{\G} \, (p_1 , \ldots , p_N) = \int d^d \, k_1 \ldots d^d \, k_L \ I_{\G}
\, (p_1 , \ldots , p_N , k_1 , \ldots , k_L) \leqno (43)
$$
\`a chaque graphe $\G$. Nous les utilisons en m\'etrique Euclidienne pour \'eviter 
les facteurs imaginaires.

\smallskip
\noindent Pour respecter les dimensions physiques des quantités impliquées 
quand on \'ecrit ces r\`egles en dimension $d$, 
il faut introduire une unit\'e de masse $\mu$ et remplacer partout la constante de
couplage par $\mu^{3-d/2} \, g$. On normalise ainsi les calculs par,
$$
U \, (\G) = g^{(2-N)} \, \mu^{-B} \, \lgl \s , U_{\G} \rgl \leqno (44)
$$
o\`u $B = B \, (d)$ est la dimension de $\lgl \s , U_{\G} \rgl$.

\smallskip

\smallskip

\noindent On étend la d\'efinition (44) aux r\'eunions disjointes de graphes 1PI
 $\G_j$ par, $$ U \, (\G = \cup \, \G_j) = \Pi \, U \, (\G_j)
\, . \leqno (45) $$ 
\smallskip

\noindent Le r\'esultat principal est alors le suivant:

\medskip

\noindent {\tencsc Théorème (\cite{ck3})} a) {\it Il existe une unique application 
m\'eromorphe $\g (z)
\in G$, $z \in \Cb$,  $z \ne D$ dont les
$\G$-coordonn\'ees sont donn\'ees par $U \, (\G)_{d=z}$.}

\smallskip

\noindent b) {\it La valeur renormalis\'ee d'une observable physique
$\Oc$ est obtenue en rempla\c{c}ant $\g \, (D)$ dans le d\'eveloppement 
perturbatif de $\Oc$ par $\g_+ \, (D)$ o\`u $$ \g \, (z) = \g_- \,
(z)^{-1} \ \g_+ \, (z) $$ est la d\'ecomposition de Birkhoff du lacet
 $\g \, (z)$ relativement à un cercle 
infinitésimal autour de
$D$.}

\section  { Le groupe de renormalisation}

\noindent  Montrons comment le groupe de renormalisation apparait très simplement 
de notre point de vue.

\noindent Comme nous l'avons vu ci-dessus, la r\'egularisation dimensionnelle 
implique le
 choix arbitraire d'une unit\'e de masse $\mu$ et l'on constate d'abord que la 
partie singulière de la d\'ecomposition de 
Riemann-Hilbert de $\g$ est en fait ind\'ependante de ce choix.
Il en r\'esulte une contrainte très forte sur cette partie singulière et le groupe 
de renormalisation s'en d\'eduit imm\'ediatement.
 Nous en d\'eduisons \'egalement une formule explicite pour l'action nue.
\medskip
\noindent On montre d'abord, en se limitant \`a la th\'eorie $\vp_6^3$ pour 
simplifier les notations,
que bien que le lacet $\g (d)$ d\'epende du choix de l'unit\'e de masse 
 $\mu$, 
$$ \mu \ra \g_{\mu} (d)
\, , \leqno (46) 
$$ 
la partie singulière $\g_{\mu^-}$ de sa d\'ecomposition de Birkhoff,
$$ \g_{\mu} (d) = \g_{\mu^-} (d)^{-1} \, \g_{\mu^+}
(d) \leqno (47) 
$$ 
est en fait ind\'ependante de $\mu$,
 $$
\frac{\partial}{\partial \mu} \, \g_{\mu^-} (d) = 0 \, . \leqno
(48) 
$$
 \noindent Cet \'enonc\'e d\'ecoule imm\'ediatement de l'analyse dimensionnelle.

\noindent De plus, par construction le groupe de Lie $G$ est muni d'un  
groupe \`a un paramètre d'automorphismes,
 $$ \t_t
\in {\rm Aut} \, G \ , \quad t \in \Rb \, , \leqno (49) $$
associ\'e \`a la graduation de l'alg\`ebre de Hopf $\Hc$ donn\'ee par le nombre de 
boucles,
 $$ L (\G) = \hbox{nombre de boucles}
\ \G \leqno (50) $$ pour tout graphe 1PI $\G$.

\noindent On a l'\'egalit\'e $$ \g_{e^t \mu} (d) = \t_{t \ve}
(\g_{\mu} (d)) \qquad \fl \, t \in \Rb \, , \ \ve = D-d \leqno (51)
$$
Il en r\'esulte que les lacets $\g_{\mu}$ associ\'es \`a la th\'eorie 
nonrenormalis\'ee
 satisfont la propri\'et\'e suivante: la partie singulière de leur d\'ecomposition 
de 
 Birkhoff est inchang\'ee  par l'op\'eration,
$$ \g (\ve) \ra \t_{t\ve} (\g (\ve)) \, , \leqno (52) $$
En d'autres termes, si l'on remplace $\g (\ve) $ par $\t_{t\ve} (\g
(\ve))$ l'on ne modifie pas la partie singulière de sa d\'ecomposition de Birkhoff.
On a pos\'e $$ \ve = D - d \in
\Cb \backslash \{ 0 \} \, . \leqno (53) $$ 
Nous donnons une caractérisation complète des lacets $\g (\ve)
\in G$ v\'erifiant cette propriet\'e. Cette caractérisation n'implique que la 
partie
singulière $\g_- (\ve)$ qui v\'erifie par hypothèse, $$ \g_-
(\ve) \, \t_{t \ve} (\g_- (\ve)^{-1}) \ \hbox{est convergent pour} \
\ve \ra 0 \, . \leqno (54) $$ 
Il est facile de voir que la limite de
(54) pour $\ve \ra 0$ d\'efinit un sous-groupe \`a un paramètre, $$ F_t \in G
\, , \ t \in \Rb \leqno (55) $$ et que le g\'en\'erateur $\b =
\left( \frac{\partial}{\partial t} \, F_t \right)_{t=0}$ de ce sous-groupe
est reli\'e au {\it r\'esidu} de $\g$ $$
\build{\rm Res}_{\ve = 0}^{} \g = - \left( \frac{\partial}{\partial
u} \, \g_- \left( \frac{1}{u} \right) \right)_{u=0} \leqno (56) $$
par l'\'equation, $$ \b = Y \, {\rm Res} \, \g \, , \leqno
(57) $$ o\`u $Y = \left( \frac{\partial}{\partial t} \, \t_t
\right)_{t=0}$ est la graduation.

\noindent Ceci est imm\'ediat mais notre r\'esultat (\cite{ck4}) donne la formule 
explicite 
(59) qui exprime $\g_- (\ve)$ en fonction de $ \b$.
Introduisons le produit semi-direct de l'alg\`ebre de Lie $G$ (des \'el\'ements 
primitifs de $\Hc^*$)
par la graduation. On a donc un \'el\'ement $Z_0$ tel que
 $$ [Z_0 , X] = Y(X) \qquad \fl \, X \in \hbox{Lie} \ G
\, . \leqno (58) $$ 
La formule pour 
 $\g_- (\ve)$ est alors
 $$ \g_- (\ve) = \lim_{t \ra \ify} e^{-t \left(
\frac{\b}{\ve} + Z_0 \right)} \, e^{t Z_0} \, . \leqno (59) $$
Les deux facteurs du terme de droite appartiennent au produit semi-direct
du groupe 
 $G$ par sa graduation, mais leur rapport (59) appartient au groupe
 $G$.

\noindent Cette formule montre que toute la structure des divergences est 
uniquement d\'etermin\'ee par
le r\'esidu et donne une forme forte des relations de t'Hooft \cite{tf}.

\section  { Le groupe G et les difféomorphismes}
  \noindent  Bien entendu, on pourrait facilement objecter aux développements précédents en 
arguant que le mystère de la renormalisation n'est pas completement éclairci car 
le groupe $G$ construit à partir des graphes de Feynman apparait également
mysterieux. Cette critique est completement levée par la merveilleuse relation,
basée sur la physique entre les algèbres de Hopf $\Hc$ des graphes de Feynman et celle, $\Hc_{cm}$
des difféomorphismes.

\noindent Nous montrons, dans le cas de masse nulle, que la formule qui donne 
la constante de couplage effective, 
$$ g_0 = \left( g + \quad \sum_{\hbox{\psfig{figure=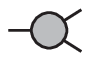}}} \quad g^{2\ell + 1} \,
\frac{\G}{S(\G)} \right) \left( 1 - \quad \sum_{\hbox{\psfig{figure=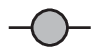}}} \quad g^{2\ell} \,
\frac{\G}{S(\G)} \right)^{-3/2} \leqno (60) $$
consid\'er\'ee comme une
s\'erie formelle dans la variable
$g$ d'\'el\'ements de l'alg\`ebre de Hopf $\Hc$, d\'efinit en fait un
 homomorphisme d'alg\`ebres de Hopf de l'alg\`ebre de Hopf 
$\Hc_{cm}$ des coordonn\'ees sur le groupe des  diff\'eomorphismes
formels de $\Cb$ tels que, $$ \vp (0) = 0 \, , \ \vp' (0) = {\rm id}
\leqno (61) $$ vers l'alg\`ebre de Hopf  $\Hc$ de la th\'eorie de masse nulle.

\noindent Il en r\'esulte en transposant, une action formelle du groupe $G$ sur la 
constante de couplage.
Nous montrons en particulier que l'image par $\rho$
de $\b = Y \, {\rm Res} \, \g$ est la fonction $\b$ de la constante de 
couplage $g$.

\noindent Nous obtenons ainsi un corollaire du th\'eor\`eme principal qui se 
formule 
sans faire intervenir ni le groupe $G$ ni l'alg\`ebre de Hopf $\Hc$.

\medskip
\noindent {\tencsc Théorème  \cite{ck4}} 
 {\it Considérons  la constante de couplage effective nonrenormalis\'ee 
 $g_{\rm eff} (\ve)$ comme une série formelle en
 $g$ et soit $g_{\rm eff}(\ve) =  g_{{\rm eff}_+} (\ve)\, (g_{{\rm 
eff}_-}(\ve))^{-1}$ 
sa décomposition de Birkhoff (opposée) dans le groupe des 
difféomorphismes formels. Alors le lacet $g_{{\rm eff}_-} (\ve)$ est la constante 
de couplage nue
et $g_{{\rm eff}_+} (0)$ la constante de couplage renormalis\'ee.}
\medskip

\noindent Comme la décomposition de Birkhoff d'un lacet à valeurs dans le groupe 
des 
difféomor-phismes (formels) est évidemment reliée à la classification des fibrés 
(non-vectoriels)
holomorphes, ce résultat suggère qu'un tel fibré ayant pour base un voisinage de 
la 
dimension $d$ de l'espace temps et pour fibre les valeurs (compléxifiées) des 
constantes de couplage
devrait donner une interprétation géométrique de l'opération de renormalisation. 
Il faut tout de même 
noter que la décomposition de Birkhoff a lieu ici relativement à un cercle 
infinitésimal autour de
$d$ et qu'il s'agit de difféomorphismes formels.

\noindent Les résultats ci-dessus montrent qu'au niveau des développements perturbatifs
le procédé de renormalisation admet une interprétation géométrique simple grâce au 
groupe $G$ et à la décomposition de Riemann-Hilbert. Le problème essentiel consiste
à passer du développement perturbatif à la théorie non-perturbative, ce qui
revient en termes de difféomorphismes à passer du developpement de Taylor
à la formule globale.

\section  { Le groupe de renormalisation et la théorie de Galois aux places 
archimédiennes}

\noindent Le problème de Riemann-Hilbert joue un rôle clef dans la théorie de Galois différentielle, il est 
donc naturel d'interpréter en termes Galoisiens l'ambiguïté que le groupe de renormalisation introduit dans les 
théories physiques.  Cette section contient l'esquisse d'une telle interprétation. Nous montrerons en 
particulier le rôle que le 
groupe de renormalisation
devrait jouer pour comprendre la composante  connexe du groupe des classe 
d'idèles de la théorie du corps de classe comme un groupe de Galois.

\noindent Commençons par une introduction très élémentaire à la théorie de Galois pour les équations 
algébriques.

\noindent Si la technique de r\'esolution des \'equations du second
degr\'e remonte \`a la plus haute Antiquit\'e (Babyloniens,
Egyptiens...), elle n'a pu \^etre \'etendue au troisi\`eme degr\'e que
bien plus tard, et ne sera publi\'ee par Girolamo (J\'er\^ome) Cardano
qu'en 1545 dans les chapitres 11 \`a 23 de son livre {\sl Ars magna sive de
regulis algebraicis}. Bien que cela n'ait pas \'et\'e reconnu avant le 
dix-huiti\`eme si\`ecle, la  clef de la r\'esolution par radicaux de l'\'equation 
g\'en\'erale 
du troisi\`eme degr\'e, $x^3+nx^2+px+q=0$, de racines $a$, $b$, $c$,
est l'existence d'une fonction rationnelle $\alpha(a,b,c)$ de $a$, $b$, $c$,
qui ne prend que deux d\'eterminations diff\'erentes sous l'action des six
 permutations de $a$, $b$, $c$.

\noindent La m\'ethode de Cardan revient \`a poser $\alpha
=\big((1/3)(a+bj+cj^2)\big)^3$ o\`u le nombre $j$ est la premi\`ere
racine cubique de l'unit\'e.
La permutation circulaire
transformant $a$ en $b$, $b$ en $c$ et $c$ en $a$ laisse
 manifestement $\alpha$  inchang\'ee et la seule autre d\'etermination
de la fonction $\alpha$ sous l'action des six permutations de $a$, $b$,  $c$,
est obtenue en transposant $b$ et $c$ par exemple,
ce qui donne 
$\beta=\big((1/3)(a+cj+bj^2)\big)^3$.

\noindent Comme l'ensemble de ces deux
nombres $\alpha$ et $\beta$ est invariant par toutes les permutations de
$a$, $b$, $c$, le polyn\^ome du second degr\'e dont $\alpha$ et
$\beta$ sont racines se calcule rationnellement en fonction des coefficients
de l'\'equation initiale $x^3+nx^2+px+q=0$
: c'est $X^2+2qX-
p^3=(X+q+s)(X+q-s)$ o\`u $s$ est l'une des racines carr\'ees de $p^3+q^2$ et 
o\`u l'on a r\'e\'ecrit l'\'equation initiale sous la 
forme \'equivalente
$x^3+3px+2q=0$ débarrassée du terme du deuxi\`eme degr\'e en
effectuant une translation convenable des racines et o\`u l'on a
introduit les coefficients 2 et 3 pour simplifier les formules.
 
\noindent Un calcul simple montre alors que chacune des racines $a$, $b$ et $c$, 
de
l'\'equation initiale s'exprime comme somme de l'une des trois racines
cubiques de $\alpha$ et de l'une des trois racines cubiques de $\beta$, ces
deux choix \'etant li\'es par le fait que leur produit doit \^etre
imp\'erativement \'egal \`a $-p$ (il n'y a donc que trois couples de choix
de ces racines \`a prendre en compte, ce qui est rassurant, \`a la place
des neuf possibilit\'es que l'on aurait pu envisager {\it a priori}).

\noindent C'est \`a l'occasion de ces formules que l'utilisation des nombres 
complexes
s'est impos\'ee. En effet, m\^eme dans le cas o\`u les trois racines sont 
r\'eelles,
il se peut que $p^3+q^2$ soit n\'egatif et que $\alpha$ et $\beta$ soient
n\'ecessairement
des nombres complexes.

\noindent Si la r\'esolution des \'equations du troisi\`eme degr\'e que nous 
venons
d'exposer a \'et\'e tr\`es longue \`a \^etre mise au point (sans doute pour
au moins l'un de ses cas particuliers entre 1500 et 1515 par Scipione del
Ferro), celle du quatri\`eme degr\'e a \'et\'e plus preste \`a la suivre
puisqu'elle figure \'egalement dans l'{\sl Ars magna}\/ (chapitre 39) o\`u
Cardano l'attribue \`a son secr\'etaire Ludovico Ferrari qui l'aurait mise
au point entre 1540 et 1545 (Ren\'e Descartes en publiera une autre en
1637). 

\noindent  Ici encore, l'on peut partir d'un polyn\^ome
débarrassé d'un coefficient, annul\'e par translation, disons $X^4+pX^2+qX+r=(X-a)(X-b)(X-c)(X-d)$.

\noindent La fonction rationnelle $\alpha(a,b,c,d)$ la plus simple\footnote{{\it Voir \cite{louis} pour 
l'ubiquité de la symétrie en question, et son rôle dans l'organisation des tournois de football}},
 ne prenant que trois 
déterminations différentes sous l'action des vingt-quatre permutations de $a$, $b$, $c$ et $d$, est 
$\alpha=ab+cd$. Les deux autres déterminations sont $\beta=ac+bd$, 
$\gamma=ad+bc$. Ce sont donc les racines d'une \'equation du troisi\`eme degr\'e
dont les coefficients s'expriment
rationnellement en fonction de $p$, $q$ et $r$. Un calcul simple montre que
 le polyn\^ome $(X-\alpha)(X-\beta)(X-
\gamma)$ est \'egal \`a $X^3-pX^2-4rX+(4pr-q^2)$. Il peut donc \^etre
d\'ecompos\'e comme on l'a vu plus haut pour en d\'eduire $\alpha$, $\beta$
et $\gamma$; en fait, il suffit m\^eme de calculer l'une seulement de ces
racines, disons $\alpha$, pour pouvoir en d\'eduire $a$, $b$, $c$ et $d$
(nous connaissons alors en effet la somme $\alpha$ et le produit $r$ des
deux nombres $ab$ et $cd$, donc ces deux nombres eux-m\^emes par une
\'equation du second degr\'e, et il ne reste plus qu'\`a exploiter les
\'egalit\'es $(a+b)+(c+d)=0$ et $ab(c+d)+cd(a+b)=-q$ pour pouvoir en
d\'eduire $a+b$ et $c+d$, donc enfin $a$, $b$, $c$ et $d$
par une autre
\'equation du second degr\'e).

\noindent C'est \`a Joseph Louis Lagrange en 1770 et 1771 (publication en 1772, 
mais
aussi, dans une moindre mesure, \`a Alexandre Vandermonde dans un m\'emoire
publi\'e en 1774 mais \'egalement r\'edig\'e vers 1770, ainsi qu'\`a Edward
Waring dans ses {\sl Meditationes algebric\ae}\/ de 1770 et \`a Francesco
Malfatti) que l'on doit la mise en lumi\`ere du r\^ole fondamental des
permutations sur les racines $a$, $b$, $c$... et sur les quantit\'es
auxiliaires $\alpha$, $\beta$..., d'ailleurs aujourd'hui justement
appel\'ees "r\'esolvantes de Lagrange".

\noindent Ces r\'esolvantes ne sont pas uniques, et par exemple le choix
$\alpha=(a+b-c-d)^2$ correspond 
\`a la solution de Descartes, mais elles fournissent la clef de toute les
r\'esolutions g\'en\'erales par radicaux. Il y en a une qui est particulièrement belle car elle est covariante
pour le groupe affine, c'est à dire vérifie l'égalité,
$$
 \alpha(\lambda \,a +z,\lambda\, b +z,\lambda\, c +z,\lambda\, d +z)=\;\lambda \,\alpha(a,b,c,d)+z
$$
 et admet donc une interprétation géométrique.
 Elle est
donnée algèbriquement par 
$$
\alpha=\frac{ad-bc}{a+d-b-c}
$$
et correspond géométriquement (figure 2) au point d'intersection des cercles circonscrits aux
triangles $ABJ$ et $JCD$ où $J$ désigne le point d'intersection des 
droites $AC$ et $BD$. 

\medskip

$$
\hbox{ 
\psfig{figure=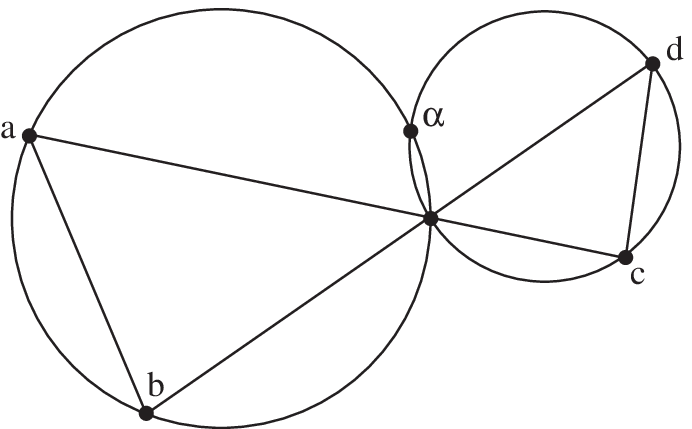} 
} 
$$

\noindent {\footnotesize Figure 2.  Le point $\alpha$ est fonction méromorphe et séparément homographique des 
quatre points A, B, C, D.}

\bigskip

\noindent J'ai rencontré récemment cette
résolvante à propos 
du problème\footnote{{\it posé par le président Chinois Jiang Zemin à
la délégation de mathématiciens venue à sa rencontre en l'an 2000.}}
 de l'étoile à cinq branches (figure 3), dont elle 
permet une résolution algébrique que je laisse à la sagacité du lecteur.

$$
\hbox{ 
\psfig{figure=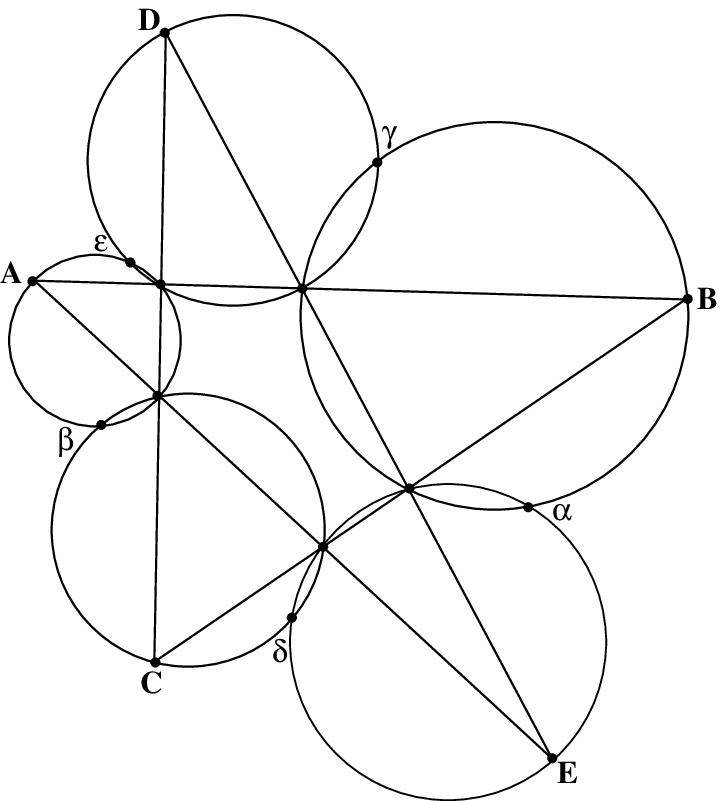} 
} 
$$

\medskip

\noindent {\footnotesize Figure 3.  On donne cinq points arbitraires  A,B,C,D,E. Montrer que les points 
d'intersection 
$ \alpha,\beta,\gamma,\delta,\epsilon$ des cercles  circonscrits aux triangles externes consécutifs de 
l'étoile sont situés sur un même cercle.}

\bigskip

\noindent L'\'etape suivante dans la théorie des équations algébriques est
\'evidemment celle du cinqui\`eme degr\'e. Descartes a certainement
essay\'e et avec lui bien des chercheurs.  Elle a toujours oppos\'e des
obstacles infranchissables, et nous savons depuis Abel et Galois, aux
alentours de 1830, pourquoi cette qu\^ete \'etait vaine.

\noindent  Descartes par exemple,
persuad\'e qu'il n'existait pas de formule analogue \`a celle de Cardano,
avait propos\'e en 1637, dans {\sl La G\'eom\'etrie}, une m\'ethode
graphique de r\'esolution gr\^ace \`a l'intersection de cercles et de
cubiques qu'il avait invent\'ees pour l'occasion. Entre 1799 et 1813 (date
de l'\'edition de ses {\sl Riflessioni intorno alla solutione delle
equazioni algebraiche generali}), Paolo Ruffini a publi\'e diverses
tentatives de d\'emonstrations, de plus en plus affin\'ees, visant \`a
\'etablir l'impossiblit\'e de r\'esoudre l'\'equation g\'en\'erale du cinqui\`eme degr\'e par radicaux. \`A toute fonction rationnelle des
racines, il a eu l'id\'ee juste d'associer le groupe des permutations de
ces racines qui la laissent invariante, mais a cru \`a tort (d'apr\`es un
rapport de Ludwig Sylow) que les radicaux intervenant dans la r\'esolution
de l'\'equation, comme les racines cubiques de $\alpha$ pour le degr\'e
trois, \'etaient n\'ecessairement des fonctions rationnelles des racines.

\noindent Il faudra attendre 1824 pour que Niels Abel justifie l'intuition de 
Ruffini
dans son {\sl M\'emoire sur les \'equations alg\'ebriques}\/ et -
apr\`es avoir cru trouver au contraire une m\'ethode de r\'esolution
g\'en\'erale - prouve  l'impossiblit\'e de r\'esoudre l'\'equation g\'en\'erale du
cinqui\`eme degr\'e par radicaux, en 1826 dans le {\sl M\'emoire
sur une classe particuli\`ere d'\'equations r\'esolubles alg\'ebriquement},
o\`u il amorce une th\'eorie g\'en\'erale qui ne s'\'epanouira que dans les 
\'ecrits
 de Galois, vers 1830. Les travaux de Galois inaugurent une \`ere nouvelle des 
math\'ematiques, o\`u les calculs font place \`a la r\'eflexion sur leur 
potentialit\'e, et 
les concepts, tels celui de groupe abstrait ou d'extension alg\'ebrique, occupent 
le devant de la sc\`ene.

\noindent L'idée lumineuse de Galois consiste d'abord à associer à une équation 
arbitraire 
un groupe de permutations qu'il définit de la manière suivante, \cite{galois}

\medskip
\noindent { \it Soit une équation donnée, dont $a$, $b$, $c$,... sont les $m$ 
racines.
Il y aura toujours un groupe de permutations des lettres $a$, $b$, $c$,...qui 
jouira de 
la propriété suivante:

1) que toute fonction des racines, invariante par les substitutions de ce groupe, 
soit rationnellement connue ;

2) réciproquement, que toute fonction des racines, déterminée
rationnellement, soit invariante par ces substitutions.}

\medskip
\noindent puis à étudier comment ce groupe "d'ambiguïté" se trouve modifié par 
l'adjonction de quantités auxiliaires
considérées comme "rationnelles".
Ainsi, dans le cas de l'équation du quatrième degré, si l'on adjoint
la quantité $\alpha$ obtenue en résolvant l'équation auxiliaire
du troisième degré, l'on réduit le groupe d'ambiguïté au sous-groupe normal formé des quatre permutations 
(a,b,c,d), (b,a,d,c),
(c,d,a,b),  (d,c,b,a). Ce groupe est le produit de deux groupes
à deux éléments et l'adjonction des solutions de deux équations du second degré suffit
alors pour éliminer totalement l'ambiguïté, c'est à dire résoudre l'équation initiale.

\noindent Si l'on désigne par $k$ le corps des "quantités rationnelles" et par
K celui engendré par $k$ et par toutes\footnote{{\it Il ne suffit pas d'adjoindre une seule de ces racines, il 
faut les adjoindre toutes}} les racines de l'équation que l'on se propose
de résoudre, le groupe de Galois, $G={\rm Gal}(K:k)$ est le groupe des 
automorphismes de $K$ qui fixent tous les éléments de $k$.

\noindent L'impossibilit\'e de r\'eduire l'\'equation du cinqui\`eme degr\'e \`a 
des \'equations 
de degr\'e inf\'erieur provient alors de la "simplicit\'e" du groupe $A_5$ des 
soixantes permutations
paires (produits d'un nombre pair de transpositions) des cinq racines $a$, $b$, 
$c$, $d$, $e$
d'une telle \'equation.
Un groupe abstrait fini est "simple" si l'on ne peut
le réduire, par un homomorphisme non trivial, à un 
groupe plus petit.
Le groupe $A_5$ est le plus petit groupe simple non commutatif
et il appara\^it très souvent en math\'ematiques.

\noindent  J'en viens maintenant au rôle que le 
groupe de renormalisation
devrait jouer pour comprendre la composante  connexe du groupe des classes 
d'idèles de la théorie du corps de classe comme un groupe de Galois.

\noindent La théorie du corps de classe et sa généralisation aux groupes de Galois non commutatifs par le 
programme de Langlands constituent l'information la plus profonde que nous ayons sur le  groupe de Galois des 
nombres algébriques.
 La source de ces théories est la loi de réciprocité quadratique qui joue un rôle central dans l'histoire
 de la théorie des nombres. Elle est démontrée en 1801 par Gauss dans ses " 
Disquisitiones " mais son énoncé était déjà connu d'Euler et de Legendre. La loi 
de réciprocité exprime,
 étant donnés deux nombres premiers $p$ et $q$, une symétrie entre
 $p$ et $q$ dans la résolution de l'équation $x^2 = p$ modulo $q$.
 Elle montre par exemple que pour savoir si l'équation $x^2 = 5$ 
admet une solution modulo un nombre premier $q$ il suffit de connaître
 la valeur de $q$ modulo $5$ ce que donne le dernier chiffre de $q $ dans son développement décimal, 
(par exemple 19 et 1999, ou 7 et 1997 donnent le même résultat)
 de sorte que les nombres premiers ainsi sélectionnés
 se répartissent en classes. Il a fallu plus d'un siècle
 pour comprendre conceptuellement la loi de réciprocité
 quadratique dont Gauss avait donné plusieurs démonstrations,
 sous la forme de la théorie du corps de classe qui permet
 de calculer à partir de classes de nombres idéaux le groupe
 de Galois de l'extension Abelienne maximale d'un corps de nombres.

\noindent La généralisation conceptuelle de la notion de corps
de nombres est celle de corps global. Un corps $k$ est {\it global} 
si c'est un sous-corps 
discret cocompact
 d'un anneau localement compact  
(non discret) semi-simple et  commutatif  $A$. (Cf. 
Iwasawa {\it Ann. of Math.} {\bf 57} (1953).) L'anneau topologique $A$ est
alors canoniquement associé à $k$ et s'appelle l'anneau des  Adèles de $k$, on a,
$$
A = \prod_{\rm res}^{} k_v \, , \leqno (62)
$$
où le produit est le produit restreint des corps locaux $k_v$ 
indéxés par les places de $k$. Les $k_v$ sont les corps localement compacts
obtenus comme complétions de $k$ de même que l'on obtient les nombres réels en complétant les rationnels.

\noindent Quand la caractéristique de $k$ est $p>1$ i.e. quand $k$ est un corps de 
fonctions sur 
 $\Fb_q$, on a
$$
k \subset k_{\rm un} \subset k_{\rm ab} \subset k_{\rm sep} \subset 
\overline k \, , \leqno (63)
$$
où $\overline k$ désigne une cloture algébrique de $k$, $k_{\rm 
sep}$ la cloture algébrique séparable, $k_{\rm ab}$ l'extension abélienne maximale
et  $k_{\rm un}$,  extension non ramifiée maximale, est obtenue en adjoignant à $k$ 
les racines de l'unité 
d'ordre premier à $p$.

\noindent On définit le groupe de Weil $W_k$ comme le sous-groupe de ${\rm Gal} 
(k_{\rm ab} : k)$ formé par les automorphismes de $(k_{\rm ab} : k)$ qui 
induisent sur $k_{\rm un}$ une puissance entière de l'automorphisme de 
"Frobenius",  
 $\theta$,
$$
\theta (\mu) = \mu^q \qquad \forall \, \mu \ \hbox{racine de l'unité d'ordre 
premier à  
} \ p \, . \leqno (64)
$$
Le résultat principal de la théorie du corps de classe global est l'existence d'un 
isomorphisme 
canonique,
$$
W_k \simeq C_k = GL_1 (A) / GL_1 (k) \, , \leqno (65)
$$
de groupes localement compacts.

\noindent Quand $k$ est de caractéristique nulle, i.e. un corps de nombres, on 
a un isomorphisme 
canonique,
$$
{\rm Gal} (k_{\rm ab} : k) \simeq C_k / D_k \, , \leqno (66)
$$
où $D_k$ désigne la composante connexe de l'élément neutre dans le groupe des classes 
d'idèles
 $C_k = GL_1 (A) / GL_1 (k)$, mais à cause des places Archimédiennes  
de $k$ l'on n'a pas d'interprétation de $C_k$ analogue au cas des corps de 
fonctions. 
Citons  A.~Weil \cite{W4}, 

\medskip

``La recherche d'une interpr\'{e}tation pour $C_k$ si $k$ est un corps de nombres, 
analogue en quelque mani\`{e}re \`{a} l'interpr\'{e}tation par un groupe de Galois 
quand $k$ est un corps de fonctions, me semble constituer l'un des 
probl\`{e}mes fondamentaux de la th\'{e}orie des nombres \`{a} l'heure actuelle~; 
il 
se peut qu'une telle interpr\'{e}tation renferme la clef de l'hypoth\`{e}se de 
Riemann~$\ldots$''.
\medskip

\noindent  Cela signifie qu'aux places Archimédiennes (i.e. aux complétions de $k$ 
qui donnent soit
les nombres réels soit les nombres complexes), il devrait y avoir un groupe 
continu de symétries
 secrètement présent.

\noindent   Mon intér\^et pour ce problème vient de mon travail
 sur la classification des facteurs qui indiquait
clairement que l'on possédait là l'analogue de la théorie de Brauer qui est l'une 
des clefs de la 
théorie du corps de classe local. 
 
\noindent Les groupes de Galois sont par construction des limites
projectives de groupes finis attachés à des extensions finies. Pour obtenir des 
groupes 
connexes il faut évidemment relaxer cette condition de finitude, qui est la même 
que 
la restriction en théorie de Brauer aux algèbres simples centrales de dimension 
finie. 
Comme ce sont les places Archimédiennes de $k$ qui sont à l'origine de la 
composante connexe $D_k$,
il est naturel de considérer la question préliminaire suivante,
\medskip

\noindent ``Existe-t-il une théorie de Brauer non-triviale d'algèbres simples 
centrales sur $\Cb$.''

\medskip

\noindent  J'ai montré dans \cite{ccc} que la classification des facteurs {\it 
approximativement finis 
} sur $\Cb$ donnait une réponse satisfaisante à cette question.
Ils sont classifiés par leur module,
$$
{\rm Mod} (M) \mathop{\subset}_{\sim} \ \Rb_+^* \, , \leqno (67)
$$
qui est un sous-groupe (virtuel) fermé de $\Rb_+^*$.

\noindent Ce groupe joue un rôle analogue dans le cas Archimédien au module des 
algèbres simples centrales
sur un corps local nonarchimédien. Dans ce dernier cas le module se définit trés 
simplement par l'action du groupe multiplicatif
d'une algèbre simple centrale sur la mesure de Haar du groupe additif. La 
définition de ${\rm Mod} (M)$ pour les 
facteurs est beaucoup plus élaborée, mais reste basée sur l'action du groupe $\Rb_+^*$ de changement d'échelle.
  
\noindent Pour poursuivre l'analogie avec la théorie de Brauer
où le lien avec le groupe de Galois s'obtient par la construction
d'algèbres simples centrales comme produits croisés d'un corps par un groupe 
d'automorphismes, le pas suivant 
consiste à trouver des exemples naturels de construction de 
facteurs comme produits croisés d'un corps  $K$, 
extension transcendante de $ \Cb$ par un groupe d'automorphismes.
Dans nos recherches sur les variétés sphériques noncommutatives
 \cite {connes:09}, avec M. Dubois-Violette, 
l'algèbre de Sklyanin (\cite{sk1})
est apparue comme solution en dimension 3, du problème de classification formulé 
dans \cite {connes:08}.
La représentation "régulière" de cette algèbre engendre 
une algèbre de von-Neumann
intégrale directe de facteurs approximativement 
finis de type ${\rm II}_1$, tous isomorphes au facteur hyperfini $R$. 
Les homomorphismes correspondants de 
l'algèbre de Sklyanin (\cite {sk1}) vers le facteur $R$ se 
factorisent miraculeusement à travers le 
produit croisé du corps  $K_q$
des fonctions elliptiques, où le 
module $q= e^{2 \pi i \tau }$ est réel,  
par l'automorphisme de translation par 
un nombre réel (mais en général irrationnel). 
 On obtient ainsi  
le facteur $R$ comme produit croisé de $K_q$ par un sous-groupe 
du groupe de Galois, en parfaite analogie avec la construction 
des algèbres simples centrales sur un corps local. 
Il reste  à obtenir une construction semblable et naturelle du facteur $R_{\infty 
}$
de type ${\rm III}_1$.

\noindent Il est sans doute prématuré d'essayer d'identifier 
le corps $K$ correspondant, qui devrait jouer le rôle de 
l'extension nonramifiée maximale $\Cb_{un}$ de $\Cb$ et être doté 
d'une action naturelle du groupe multiplicatif $\Rb_+^*$.

\noindent Le rôle du corps $K$ en physique des hautes énergies devrait être relié à l'observation 
suivante concernant les "constantes"
qui interviennent en théorie des champs.  En fait les calculs des physiciens regorgent d'exemples de 
"constantes"
telles les constantes de couplage $g$ des interactions (électromagnétiques, 
faibles et fortes)
 qui n'ont de "constantes"
que le nom. Elles dépendent en réalité du niveau d'énergie $\mu$ 
auquel les expériences sont réalisées et sont donc des fonctions $g( \mu)$. Ainsi les physiciens des hautes 
énergies
étendent implicitement le "corps des constantes" avec lequel
 ils travaillent, passant du corps $\Cb$ des scalaires 
à un corps de fonctions $g( \mu)$.
 Le générateur du groupe de renormalisation est simplement $ \mu \partial 
/ \partial \mu$.

\noindent L'on peut mettre l'exemple plus simple du corps $K_q$ des fonctions elliptiques 
sous la même forme en passant aux fonctions loxodromiques, c'est à dire 
en posant $\mu = e^{2 \pi i z }$ de sorte que la première périodicité (en $z 
\rightarrow z+1$) est automatique
alors que la deuxième s'écrit $g( q \, \mu)=g( \mu) $. Le groupe des automorphismes de la 
courbe elliptique est alors lui aussi engendré par $ \mu \partial / \partial \mu$.

\noindent Les points fixes du groupe 
de renormalisation sont les scalaires ordinaires, mais il se pourrait que la 
physique quantique
conspire pour nous empêcher d'espérer une théorie qui englobe
toute la physique des particules et soit construite comme point 
fixe du groupe 
de renormalisation. Les interactions fortes sont asymptotiquement libres et l'on peut les analyser
 à très hautes énergies en utilisant les points fixes du groupe
de renormalisation, mais la présence du secteur électrodynamique montre
 qu'il est vain de vouloir s'en tenir à de tels points fixes pour décrire une 
théorie qui incorpore 
l'ensemble des forces observées. Le problème est le même dans le domaine
infrarouge où les rôles des interactions fortes et électrodynamiques sont 
inversés.

\noindent Il est bien connu des physiciens que le groupe de renormalisation joue 
le rôle d'un groupe 
d'ambiguïté, l'on ne peut distinguer entre elles deux théories physiques qui appartiennent à 
la 
même orbite de ce groupe, ce qui nous ramène à  Galois dont
la "théorie de l'ambiguïté" allait bien au delà des équations algébriques.  

 \noindent Citons \'Emile Picard (voir \cite{ramis}) qui dans 
 sa préface aux \oe uvres complètes d'Evariste  Galois 
écrivait,

\noindent "Il aurait édifié dans ses parties essentielles, la théorie des fonctions 
algébriques d'une variable telle que nous la connaissons aujourd'hui. Les 
méditations de Galois portèrent
encore plus loin; il termine sa lettre en parlant de l'application à l'analyse 
transcendante de la théorie de
l'ambiguïté. On devine à peu près ce qu'il entend par là, et sur ce terrain qui, 
comme il le dit est immense, il reste encore aujourd'hui bien
des découvertes à faire".

\end{document}